\RequirePackage[2020-02-02]{latexrelease}
\documentclass[11pt,oneside]{amsart}
\usepackage[left=1in, right=1in]{geometry}
\usepackage{amsmath,amsthm,amsfonts,amssymb,amscd,mathabx}

\makeatletter
\let\save@mathaccent\mathaccent
\newcommand*\if@single[3]{%
	\setbox0\hbox{${\mathaccent"0362{#1}}^H$}%
	\setbox2\hbox{${\mathaccent"0362{\kern0pt#1}}^H$}%
	\ifdim\ht0=\ht2 #3\else #2\fi
}
\newcommand*\rel@kern[1]{\kern#1\dimexpr\macc@kerna}
\renewcommand*\widebar[1]{\@ifnextchar^{{\wide@bar{#1}{0}}}{\wide@bar{#1}{1}}}
\newcommand*\wide@bar[2]{\if@single{#1}{\wide@bar@{#1}{#2}{1}}{\wide@bar@{#1}{#2}{2}}}
\newcommand*\wide@bar@[3]{%
	\begingroup
	\def\mathaccent##1##2{%
		\let\mathaccent\save@mathaccent
		\if#32 \let\macc@nucleus\first@char \fi
		\setbox\z@\hbox{$\macc@style{\macc@nucleus}_{}$}%
		\setbox\tw@\hbox{$\macc@style{\macc@nucleus}{}_{}$}%
		\dimen@\wd\tw@
		\advance\dimen@-\wd\z@
		\divide\dimen@ 3
		\@tempdima\wd\tw@
		\advance\@tempdima-\scriptspace
		\divide\@tempdima 10
		\advance\dimen@-\@tempdima
		\ifdim\dimen@>\z@ \dimen@0pt\fi
		\rel@kern{0.6}\kern-\dimen@
		\if#31
		\overline{\rel@kern{-0.6}\kern\dimen@\macc@nucleus\rel@kern{0.4}\kern\dimen@}%
		\advance\dimen@0.4\dimexpr\macc@kerna
		\let\final@kern#2%
		\ifdim\dimen@<\z@ \let\final@kern1\fi
		\if\final@kern1 \kern-\dimen@\fi
		\else
		\overline{\rel@kern{-0.6}\kern\dimen@#1}%
		\fi
	}%
	\macc@depth\@ne
	\let\math@bgroup\@empty \let\math@egroup\macc@set@skewchar
	\mathsurround\z@ \frozen@everymath{\mathgroup\macc@group\relax}%
	\macc@set@skewchar\relax
	\let\mathaccentV\macc@nested@a
	\if#31
	\macc@nested@a\relax111{#1}%
	\else
	\def\gobble@till@marker##1\endmarker{}%
	\futurelet\first@char\gobble@till@marker#1\endmarker
	\ifcat\noexpand\first@char A\else
	\def\first@char{}%
	\fi
	\macc@nested@a\relax111{\first@char}%
	\fi
	\endgroup
}
\makeatother

\makeatletter
\renewcommand*\env@matrix[1][*\c@MaxMatrixCols c]{%
	\hskip -\arraycolsep
	\let\@ifnextchar\new@ifnextchar
	\array{#1}}
\def\thm@space@setup{%
	\thm@preskip=\parskip \thm@postskip=0pt
}
\makeatother

\makeatletter
\@namedef{subjclassname@2020}{\textup{2020} Mathematics Subject Classification}
\makeatother

\usepackage{lastpage}
\usepackage{enumerate}
\usepackage{fancyhdr}
\usepackage{mathrsfs}
\usepackage{xcolor}
\usepackage{graphicx}
\graphicspath{ {./Diagrams/} } 
\usepackage{listings}
\usepackage{euscript}
\usepackage{amsthm}
\usepackage{xifthen}
\usepackage{verbatim}
\usepackage{hyperref}
\usepackage{cleveref}
\usepackage{lipsum}
\usepackage{float}
\usepackage{subfigure}
\usepackage{tikz-cd}

\usepackage{xwatermark}
\usepackage{arcs}
\usepackage{etoolbox}

\makeatletter
\providecommand\@gobblethree[3]{}
\patchcmd{\over@under@arc}
{\@gobbletwo}
{\@gobblethree}
{}{}
\makeatother

\newcommand{\overarcmath}[1]{\overarc{$#1$}}

\usepackage[utf8]{inputenc}
\usepackage[T1]{fontenc}

\usepackage{array}
\newcolumntype{M}[1]{>{\centering\arraybackslash}m{#1}}

\hypersetup{%
	colorlinks=true,
	linkcolor=blue,
	citecolor=blue,
	linkbordercolor={0 0 1}
}

\lstdefinestyle{Python}{
	language        = Python,
	frame           = lines, 
	basicstyle      = \footnotesize,
	keywordstyle    = \color{blue},
	stringstyle     = \color{green},
	commentstyle    = \color{red}\ttfamily
}

\newcommand{\mc}[1]{\mathcal{#1}}
\newcommand{\mf}[1]{\mathfrak{#1}}
\newcommand{\Z}{\mathbb{Z}}

\newcommand{\boundary}   {{\ensuremath \partial}}
\newcommand{\normalclosure}[1]{\langle\!\langle #1\rangle\!\rangle}
\newcommand{\Aut}{\operatorname{Aut}}
\newcommand{\p}{\textup{\textsf{p}}}

\newcommand{\bcup}{\bigcup}

\newcommand{\sseq}{\subseteq}

\theoremstyle{definition}
\newtheorem{theorem}{Theorem}[section]
\theoremstyle{definition}
\newtheorem{lem}[theorem]{Lemma}
\theoremstyle{definition}
\newtheorem{conj}[theorem]{Conjecture}
\theoremstyle{definition}

\theoremstyle{definition}
\newtheorem{prop}[theorem]{Proposition}
\theoremstyle{definition}
\newtheorem{cor}[theorem]{Corollary}
\theoremstyle{definition}
\newtheorem{examp}[theorem]{Example}
\theoremstyle{definition}
\newtheorem{rem}[theorem]{Remark}
\theoremstyle{definition}
\newtheorem{defin}[theorem]{Definition}
\theoremstyle{definition}
\newtheorem{construction}[theorem]{Construction}
\theoremstyle{definition}

\newcommand{\bkt}[1]{\left(#1\right)} 

\author[H.~K.~Chong]{Hip Kuen Chong}
\email{chonghk1997@gmail.com}
\author[D.~T.~Wise]{Daniel T. Wise}
\email{wise@math.mcgill.ca}
\date{\today}
\thanks{Research supported by NSERC}

\title{Embedding Partial HNN Extensions In Ascending HNN Extensions}

\begin{document}

\begin{abstract}
	We show that any partial ascending HNN extension of a free group embeds in an actual ascending HNN extension of a free group.  Moreover, we can ensure that it embeds as the parabolic subgroup of a relatively hyperbolic group.
\end{abstract}

\maketitle

\section{Introduction}

Feighn and Handel \cite{MR1709311} proved that any finitely generated group $G$ that is an ascending HNN extensions of a free group is finitely presented.  Their proof actually shows that $G$ splits as a partial ascending HNN extension.  Let $\Phi\colon F\to F$ be a monomorphism of the free group $F = \langle a_i\colon i\in I\rangle$.  The associated \emph{ascending HNN extension} is $G = F*_\Phi= \langle t, a_i\colon i\in I\mid a_i^t = \Phi(a_i)\colon i\in I\rangle$.  A \emph{partial ascending HNN extension} is a group of the form $G' = F*_\Psi= \langle t, a_i\colon i\in I\mid a_i^t = \Psi(a_i)\colon i\in J\rangle$ where $J\sseq I$ and $\Psi\colon F'\to F$ is a monomorphism where $F' = \langle a_i\colon i\in J\rangle$. 


The motivation of this paper is to prove a converse:
\begin{theorem}\label{thm:mainHNN}
	Let $G'$ be a partial ascending HNN extension of a f.g.\ free group.  Then there is an injection $G'\sseq G$ where $G$ is an ascending HNN extension of a f.g.\ free group.  Moreover, we can ensure that $G'\hookrightarrow G$ arises as $F*_\varphi\hookrightarrow F''*_\Phi$ where $\Phi$ extends $\varphi$ as follows:
	\[ \begin{tikzcd}
		F' \arrow{r}{\varphi} \arrow[swap]{d}{} & F \arrow{d}{}\\%
		F''\arrow{r}{\Phi}& F''
	\end{tikzcd}
	\]
\end{theorem}

Our result raises the following: 
\begin{conj}
	Every f.g.\ free-by-cyclic group is a subgroup of a (f.g.\ free)-by-cyclic group.  Moreover, there is an extension $\varphi\to \Phi$ as in the above diagram.
\end{conj}
This corresponds to the case where both $F'\sseq F$ and $\Phi(F')\sseq F$ are free factors.

The method we employ to prove Theorem~\labelcref{thm:mainHNN} is of independent interest.  

\begin{prop}\label{prop:sufficientinjectionCopy}
	Let $Y\sseq X$ be a subcomplex.  Suppose  $X/Y$ is combinatorially reducible and $X\to X/Y$ has liftable cancellable pairs.
	Then $\pi_1 Y\to \pi_1 X$ is injective.
\end{prop}
Proposition~\labelcref{prop:sufficientinjectionCopy}, which is proven as Proposition~\labelcref{prop:sufficientinjection}, is simple and definitional, but highly applicable.  
Proposition~\labelcref{prop:sufficientinjectionCopy} is proven by considering a disk diagram $D\to X$,
observing the quotient disk diagram $\widebar D\to X/Y$ must have a cancellable pair, and lifting this to a cancellable pair in $D\to X$.  

\begin{prop}\label{prop:sufficientrelativehypCopy}
	Let $Y\sseq X$ be a subcomplex of a compact $2$-complex $X$.  
	Suppose $X/Y$ is combinatorially reducible and $X\to X/Y$ has liftable cancellable pairs. 
	Suppose there is $K>0$ such that for each reduced disk diagram $D\to X$, the induced diagram $\widebar D\to X/Y$ satisfies $\operatorname{Area}(\widebar{D})\leq K\cdot |\boundary_\p \widebar{D}|$.  
	Then $\pi_1 X$ is hyperbolic relative to $\pi_1 Y$.  
\end{prop}
Proposition~\labelcref{prop:sufficientrelativehypCopy} is proven in Proposition~\labelcref{prop:sufficientrelativehyp} and Corollary~\labelcref{cor:sufficientrelativehypcor} with the aid of Osin's relative hyperbolicity criterion. We consider an associated disk diagram $\overarcmath{D}$ to Osin's complex $\overarcmath{X}$ and relate the numbers of cells in $\overarcmath{D}$ and the quotient diagram $\widebar D$. In particular, we 
show $\operatorname{Area}_{X-Y}(\overarcmath{D}) = \operatorname{Area}(\widebar{D})$, and show that the number of remaining $2$-cells is at most proportional to $\operatorname{Area}_{X-Y}(\overarcmath{D})$.

Theorem~\labelcref{thm:mainHNN} is proven by letting $Y$ be the mapping torus of $\varphi$, and then adding two generators to the free group and extending $\varphi$ to those new generators.  The resulting mapping cylinder is $X$.  We extend $\varphi$ by adding long relators to ensure the small-cancellation property for $X/Y$.  This is described in Construction~\labelcref{construction:halfclean_to_ascending} and $\pi_1 Y\leq \pi_1 X$ is proved in Proposition~\labelcref{prop:groupembedding}.  The extension $Y\hookrightarrow X$ is illustrated by the following presentations where the new cells are indicated in bold: 
\begin{align*}
	&\langle a,b,c,t\mid a^t = (abc)^8,\ b^t =(ac)^9b\rangle\\
	&\langle a,b,c,\mathbf{x,y},t\mid a^t = (abc)^8,\ b^t=(ac)^9b,\ \mathbf{c^t = (xy)^{100},\  x^t=x(xxy)^{100},\ y^t = y(xxxy)^{100}}\rangle
\end{align*}
We end Section~\labelcref{section:SmallCancellation} with an amusing application: if one assumes moreover that $\varphi$ is fully irreducible, then we may ensure its extension $\Phi$ is also fully irreducible.  See Definition~\labelcref{defin:irreducibleHNN} and Proposition~\labelcref{prop:IrredEmbedding}.


\section{Disk Diagram and Diagrammatic Reducibility}

Let $X$ be a combinatorial $2$-complex.  Denote $X^k$ as the $k$-dimensional skeleton of $X$. 

\begin{defin}
	The \emph{boundary path} of a $2$-cell $R$ is denoted by $\boundary_\p R$.  The boundary path is an \emph{$n$-th power} if $\boundary_\p R = w^n$ for some cycle $w\to X^1$.  
	The boundary path has \emph{exponent $n$} if it is an $n$-th power but not an $m$-th power for $m>n$.
\end{defin}

\begin{defin}
	Let $D\to X$ be a combinatorial map.  A \emph{cancellable pair} of $2$-cells $R_1, R_2$ along a $1$-cell $e$ consists of $2$-cells of $D$ where $\boundary_\p R_1 = eP_1$ and $\boundary_\p R_2 = eP_2$, and the compositions $eP_1 \to X$ and $eP_2 \to X$ are the same closed path.  
	
	A \emph{disk diagram} $D$ is a compact contractible $2$-complex with a chosen spherical embedding $D\hookrightarrow S^2$.  The \emph{boundary cycle} $\boundary_\p D$ of $D$ is the cycle equals $\boundary_\p C_\infty$, where $C_\infty$ is the open $2$-cell consisting of $S^2-D$.  
	
	A map $D\to X$ is \emph{reduced} if it has no cancellable pair.  We are especially interested in reduced maps where $D$ is a disk diagram or $D\cong S^2$ is a $2$-sphere.  We often refer to a combinatorial map $D\to X$ as a \emph{disk diagram $D$ in $X$} when $D$ is a disk diagram.  We similarly define a \emph{spherical diagram in $X$}.  The motivation for reduced disk diagrams is that if $D$ has a cancellable pair, then one can remove the open $2$-cells $R_1,R_2$ and an open arc containing $e$ between them, and then obtain a smaller spherical or disk diagram with the same boundary path by gluing.  
\end{defin}



\begin{defin} 
	A $2$-complex $X$ is \emph{combinatorially reducible} if each combinatorial map $S^2\to X$ from a $2$-sphere has a cancellable pair.  Equivalently, there is no reduced spherical diagram $S^2\to X$.  
\end{defin}

Some examples of combinatorially reducible $2$-complexes are staggered $2$-complexes, complexes satisfying the Gersten-Pride weight-test condition, and complexes satisfying the classical small-cancellation conditions (see Proposition~\labelcref{prop:C7toCombinReducible}).

\section{Subcomplex $\pi_1$-Injection}\label{sec:pi1injection}


In this section, we state a criterion for $\pi_1$-injectivity of a subcomplex $Y\sseq X$.  Consider the quotient map $\rho\colon X\to \widebar{X}$, where $\widebar{X} = X/Y$.  For a $2$-cell $R$ in $X$, let $\widebar{R}$ denote its image in $\widebar{X}$. 


\begin{rem}
	It is possible that $\boundary_\p \widebar{R}$ has backtracks even when $\boundary_\p R$ is an immersion.  
\end{rem}

\begin{construction}
	Given a disk diagram $\phi\colon D\to X$, 
	
	
	the \emph{projected diagram} $\bar \phi\colon \widebar{D}\to \widebar{X}$ is constructed by first taking $\widebar D=D\to X$, then quotient each component of the preimage of the basepoint in the map $X\to \widebar{X}$ to a point.  
	
	Let $\widebar R$ be a $2$-cell in $\widebar D$.  Let the $2$-cell $R\sseq D$ be the preimage of $\widebar R$ under the quotient $D\to \widebar D$.  We set $\bar\phi(\widebar R) = \overline{\phi(R)}$, where $\overline{\phi(R)}$ is the image of $\phi(R)$ under $X\to \widebar{X}$.
\end{construction}

\begin{defin}

	$X\to \widebar{X}$ has \emph{liftable cancellable pairs} if for each disk diagram $D\to X$ projecting to $\widebar{D}\to \widebar{X}$, if $\widebar{R}_1, \widebar{R}_2$ is a cancellable pair in $\widebar{D}$ then $R_1, R_2$ is a cancellable pair in $D$.  (It suffices to consider diagrams with two $2$-cells.) 
\end{defin}

\begin{lem}\label{lem:EquivalenceLiftableCancellablePairs}
	The following statements are equivalent:
	\begin{enumerate}[(i)]
		\item $X\to \widebar X$ has liftable cancellable pairs.
		\item For every reduced disk diagram $D\to X$, the projected diagram $\widebar D \to \widebar X$ is also reduced.
	\end{enumerate}
\end{lem}
\begin{proof}
	((i) $\implies$ (ii)) Suppose $D\to X$ is reduced and assume $\widebar D \to \widebar X$ is not reduced.  Then there is a cancellable pair $\widebar R_1, \widebar R_2$ in $\widebar D$.  Their preimages $R_1, R_2$ form a cancellable pair in $D$ by the lifting property, contradicting the reducibility of $D$.
	
	((ii) $\implies$ (i)) Suppose $\widebar R_1, \widebar R_2$ is a cancellable pair in $\widebar D$.  Then their preimages $R_1,R_2$ share at least an edge with each other.  Assume $R_1,R_2$ do not form a cancellable pair.  Let $D'\to X$ to be the disk diagram containing only $R_1, R_2$.  Then $D'$ is reduced, and hence $\widebar D' \to \widebar X$ is also reduced, leading to a contradiction.
\end{proof}

\begin{prop}\label{prop:sufficientinjection}
	Let $Y\sseq X$ be a subcomplex.  Suppose  $\widebar X$ is combinatorially reducible and $X\to \widebar X$ has liftable cancellable pairs.
	Then $\pi_1 Y\to \pi_1 X$ is injective.
\end{prop}
\begin{proof}
	Let $\gamma\to Y$ be a closed combinatorial path. We show that if $\gamma\to X$ is nullhomotopic, then $\gamma\to Y$ is nullhomotopic.  
	
	
	Suppose $\gamma = \boundary_\p D$ for some reduced disk diagram $D\to X$ (existence by \cite{VanKampen33}).  
	Assume $\gamma\to Y$ is essential, then $\widebar D$ contains at least one $2$-cell $\alpha$.  
	
	
	
	
	
	
	
	
	We claim $\widebar D$ is homeomorphic to a \emph{singular surface} obtained by gluing surfaces together along vertices and isolated edges.  Since $\boundary_\p D$ maps to $Y$, we have $\overline{\boundary_\p D}$ is a point.  The link of each vertex of $\widebar D$ is a disjoint union of circles and points.  Indeed, the preimage of a $0$-cell $v$ in $\widebar D$ is a connected subcomplex $Z\sseq D$ mapping to $Y$.  Each component of the link of $v$ corresponds to a cycle or (possibly trivial) arc in the boundary of a regular neighbourhood of $Z$.   
	
	
	
	
	We now show $\widebar{D}$ is simply connected.  
	Since we may quotient components one at a time and apply induction, it suffices to show that $H_1(D/C)=0$ for each component $C$ of the preimage of $Y$.  By the relative Mayer-Vietoris sequence $H_1(D)\to H_1(D,C)\to H_0(C)\to H_0(D)\to H_0(D,C)$, since $H_1(D) = 0 = H_0(D,C)$ and $H_0(C) = \Z = H_0(D)$, we have $H_1(D,C) = 0$.  Since $C\sseq D$ is a subcomplex, $H_1(D/C) = H_1(D,C) = 0$.  
	
	Each surface of $\widebar D$ is a sphere since it  homologically injects.  It follows that $\widebar D$ is a tree of spheres.  Consider a sphere in $\widebar D$ containing $\alpha$.  Since $\widebar{X}$ is combinatorially reducible, there is a cancellable pair $\widebar{R}_1,\widebar{R}_2\sseq \widebar{D}$.  
	%
	%
	Since $X \to \widebar{X}$ has liftable cancellable pairs, the pair of preimages $R_1,R_2\sseq D$ is cancellable.  This contradicts that $D$ is reduced.   
	%
	Hence, $Y\to X$ is $\pi_1$-injective.
\end{proof}


\section{Small Cancellation}\label{section:SmallCancellation}

\subsection{Background}

\begin{defin}
	In a group presentation $\langle\,S\,|\, R\,\rangle$, a \emph{piece} is a common cyclic subword $w$ appearing as two different ways in $R$.  Note that for a relator $r=q^n$, subwords that differ by a $\Z_n$-action are regarded as appearing in the same way.  The presentation is \emph{C(7)} if all $r\in R$ cannot be written as a concatenation of $6$ pieces.
\end{defin}

\begin{defin}
	The \emph{presentation complex} $X$ associated to a group presentation $\langle\,S\,|\, R\,\rangle$ is constructed by:
	\begin{itemize}
		\item a $0$-cell $v$;
		\item a $1$-cell $e_s$ for each $s\in S$.  We label each $e_s$ by $s$, and  attach each $e_s$ to $v$ on both sides;
		\item a $2$-cell $f_r$ for each $r\in R$, with the attaching map corresponding to the word $r$.
	\end{itemize}
\end{defin}

\begin{defin}
	\emph{Length} of a combinatorial path $p\to X$ is the number of $1$-cell in $p$, denoted by $|p|$.  The \emph{area of a combinatorial diagram $D\to X$} is the number of $2$-cell in $D$, denoted by \emph{$\operatorname{Area}(D)$}.  
	
	A finite presentation \emph{admits a linear isoperimetric function} if there is a constant $K>0$ such that for all reduced disk diagram $D$ in the associated presentation complex, $\operatorname{Area}(D) \leq K\cdot |\boundary_\p D|$.  
\end{defin}

\begin{prop}\label{prop:C7LinearIsoFunc}
	Let $G = \langle\,S\,|\, R\,\rangle$ be a finite $C(7)$ presentation.  Then it admits a linear isoperimetric function.  
\end{prop}
\begin{proof}
	For any word $w$ representing $1_G$, there is a reduced disk diagram $D$ such that $\boundary_\p D = w$.  By \cite[Cor 2.4]{Lyndon_1966}, there is a shell $R\sseq D$ such that $\boundary_\p R\cap \boundary_\p D$ has more pieces than $\boundary_\p R - \boundary_\p D$.  Hence $\operatorname{Area}(D)$ is upper bounded by the number of pieces of $w$ by removing shells in succession.
\end{proof}

\begin{prop}\label{prop:C7toCombinReducible}
	Every spherical diagram $D\to X$ to a $C(7)$ presentation complex contains a cancellable pair.  In other words, every $C(7)$ presentation complex is combinatorially reducible.
\end{prop}
\begin{proof}
	Let $D\to X$ be a reduced spherical diagram.  Let $R\sseq D\to X$ be a $2$-cell.  Since $X$ is a presentation complex and $D$ is spherical, $R\neq D$.  Hence, $D-R\to X$ is a reduced disk diagram, hence there is a shell $R'$ in $D-R$ such that $R'\cap \boundary_\p (D-R)$ is a concatenation of at least 2 pieces \cite[Cor 2.4]{Lyndon_1966}.  This is not possible since $\boundary_\p (D-R) = \boundary_\p R$ and $R$ is a single $2$-cell.
\end{proof}



\subsection{Ensuring Injectivity via Small Cancellation}

Let $Y\sseq X$ be a combinatorial subcomplex.  Let $\widebar X = X/Y$ as the beginning of Section~\labelcref{sec:pi1injection} and denote the image of a $2$-cell $R$ by $\widebar R$.

\begin{defin}
	$X$ has \emph{no extra powers relative to $Y$} if for each $2$-cell $R \not\sseq Y$, if $\boundary_\p \widebar{R}$ has exponent $n$ for some $n$, then $\boundary_\p R$ also has exponent $n$.
\end{defin}

\begin{defin}
	$X$ has \emph{no duplicates relative to $Y$} if for any $2$-cells $R_1, R_2\not\sseq Y$, if $\boundary_\p \widebar{R}_1=\boundary_\p \widebar{R}_2$ then $\boundary_\p R_1=\boundary_\p R_2$.
\end{defin}

\begin{examp}
	Let $X_1$ be the presentation complex of $\langle a, b, c \mid bcabcbc\rangle$.  Let $X_2$ be the presentation complex of $\langle a, b, c \mid abc, abcc\rangle$.  Let $Y_a$ be a subcomplex of $X_i$ corresponding to $\langle a\rangle$, and let $Y_c$ be a subcomplex of $X_i$ corresponding to $\langle c\rangle$, for $i=1, 2$.  
	
	$X_1$ has extra powers relative to $Y_a$ since $X_1/Y_a$ is the presentation complex of $\langle b, c\mid bcbcbc\rangle$, which $bcbcbc$ has exponent $3$ instead of $bcabcbc$ having exponent $1$.  On the other hand, $X_1$ has no extra powers relative to $Y_c$ since $X_1/Y_c$ is the presentation complex of $\langle a, b\mid babb\rangle$, which $babb$ also has exponent $1$.  
	
	$X_2$ has duplicates relative to $Y_c$ since $X_2/Y_c$ is the presentation complex of $\langle a, b\mid ab, ab\rangle$.  On the other hand, $X_2$ has no duplicates relative to $Y_a$ since $X_2/Y_a$ is the presentation complex of $\langle b, c\mid bc, bcc\rangle$.
\end{examp}

\begin{lem}\label{lem:LiftableCancellable}
	Suppose $X$ has no extra powers and no duplicates relative to $Y$.  Then $X\to \widebar{X}$ has liftable cancellable pairs.
\end{lem}
\begin{proof}
	Let $\widebar{R}_1, \widebar{R}_2$ be a cancellable pair intersecting at $1$-cell $\bar e$ in a disk diagram $\widebar{D}$ in $\widebar{X}$.  Then there is an isomorphism $\bar f\colon \boundary_\p\widebar{R}_1\to \boundary_\p\widebar{R}_2$ with $\bar f(\bar e)\in \Aut(\boundary_\p\widebar{R}_2)\bar e$.
	Since $X$ has no duplicates relative to $Y$, the isomorphism $\bar f$ lifts to $f\colon \boundary_\p R_1\to \boundary_\p R_2$. 
	
	Let $e$ be the preimage of $\bar e$.  Since $\Aut(\boundary_\p R_2)\cong \Aut(\boundary_\p \widebar{R}_2)$ by no extra relative powers, $f(e)\in \Aut(\boundary_\p R_2)e$.  Hence, $g\circ f$ fixes $e$ for some $g\in \Aut(\boundary_\p R_2)$ and so $R_1$ and $R_2$ form a cancellable pair along $e$.  

\end{proof}
\begin{cor}\label{cor:Injection}
	Suppose $X$ has no extra powers and duplicates relative to $Y$ and the $2$-skeleton of the quotient $\widebar X$ is a $C(7)$ presentation complex.  Then $\pi_1 Y\to \pi_1 X$ is injective.
\end{cor}
\begin{proof}
	By Proposition~\labelcref{prop:C7toCombinReducible}, $2$-skeleton of $\widebar X$ is $C(7)$ implies $\widebar X$ is combinatorially reducible since the definition of combinatorially reducible only depends on $2$-skeleton.  By Lemma~\labelcref{lem:LiftableCancellable}, $X\to \widebar X$ has liftable cancellable pairs.  By Proposition~\labelcref{prop:sufficientinjection}, $\pi_1 Y\to \pi_1 X$ is injective.
\end{proof}

\subsection{Partial Ascending HNN Extension embeds in Ascending HNN Extension}


We will show that any partial ascending HNN extension of a finitely generated free group is isomorphic to a subgroup of an ascending HNN extension of a finitely generated free group by applying Corollary~\labelcref{cor:Injection}.

\begin{defin}
	Let $G = F *_{F_1^t = F_2}$ be the HNN extension of a group $F$ associated with an isomorphism $\phi\colon F_1\to F_2$ between subgroups of $F$.  The HNN extension is \emph{ascending} if $F_1 = F$.  The HNN extension is \emph{partial ascending} if $F_1$ is a free factor of $F$.  
\end{defin}

\begin{construction}\label{construction:halfclean_to_ascending}
	Given $H$ a partial ascending HNN extension of a finitely generated free group $F$, we will construct an ascending HNN extension $G$ with $H\leq G$.  The subgroup relationship is proved in Proposition~\labelcref{prop:groupembedding}.
	
	$H$ can be presented as:
	\begin{equation}\label{eq:HalfCleanHNN_GroupH}
		H\ =\ \left\langle a_i, b_j,t\colon i\in I, j\in J \ \middle|\  ta_it^{-1} = A_i\colon i\in I\right\rangle
	\end{equation}
	where $F = \langle a_i, b_j \colon i\in I, j\in J\rangle$ and $F_1 = \langle a_i\colon i\in I\rangle$ and $A_\ell$ is a reduced word in $\{a_i^{\pm 1}, b_j^{\pm 1}\colon i\in I, j\in J\}$ for each $\ell\in I$.  
	
	Let $\{B_j, C_1,C_2\colon j\in J\}$ be a set of words in $\{c_1,c_2\}$ such that $W=\{B_j, c_1^{-1}C_1,c_2^{-1}C_2\colon j\in J\}$ satisfies $C(7)$ in the sense that $\langle c_1, c_2 | W\rangle$ is $C(7)$ and no word in $W$ is a proper power. 
	We construct $G$ to be an ascending HNN extension with presentation
	\begin{equation}\label{eq:HalfCleanHNN_GroupG}
		G\ =\ \left\langle a_i, b_j, c_1, c_2,t \colon i\in I, j\in J\ \middle|\  ta_it^{-1} = A_i,\ tb_jt^{-1} = B_j,\ tc_1t^{-1} = C_1,\ tc_2t^{-1} = C_2\colon i\in I, j\in J\right\rangle
	\end{equation}
\end{construction}

\begin{prop}\label{prop:groupembedding}
	With $H$ and $G$ as in Construction~\labelcref{construction:halfclean_to_ascending}, there is an injective homomorphism from $H$ to $G$.
\end{prop}
\begin{proof}
	Let $X$ and $Y$ be the presentation complexes of $G$ and $H$ respectively.  We regard $Y$ as a subcomplex of $X$.  It suffices to show that $\pi_1 Y\to \pi_1 X$ is injective.  Recall $\widebar X = X/Y$.  
	
	$\widebar{X}^1 = \{c_1,c_2\}$ and $2$-cells in $\widebar{X}$ have attaching maps labelled by $W$. Hence, $\widebar{X}$ is $C(7)$.
	
	Since words in $W$ have exponent $1$, $X$ has no extra power relative to $Y$.  
	
	Since words in $W$ are distinct, $X$ has no duplicates relative to $Y$.  
	
	Therefore, $H = \pi_1 Y\to \pi_1 X = G$ is injective by Corollary~\labelcref{cor:Injection}.
\end{proof}

\subsection{Partial Ascending Fully Irreducible HNN Extension embeds in Fully Irreducible Ascending HNN Extension}

In this subsection, we show that any partial ascending fully irreducible HNN extension of a finitely generated free group is isomorphic to a subgroup of a fully irreducible ascending HNN extension of a finitely generated free group by applying Corollary~\labelcref{cor:Injection}.  We refer to \cite{abdenbi_2023} for more about irreducible partial endomorphisms and their relationship with other combinatorial group theory notions.

\begin{defin}\label{defin:irreducibleHNN}
	Let $G = F *_{F_1^t = F_2}$ be the HNN extension associated with an isomorphism $\phi\colon F_1\to F_2$ between subgroups of $F$. 
	The isomorphism $\phi$ is \emph{weakly reducible} if there are $p>1$, a proper nontrivial free factor $A\sseq F$ and $g\in F$ such that $\phi^p(A)$ is well-defined and $\phi^p(A)\sseq gAg^{-1}$.  Here $A$ is a \emph{weakly invariant free factor} of $\phi$.
	The isomorphism $\phi$ is \emph{fully irreducible} if $\phi$ is not weakly reducible.  $G$ is \emph{fully irreducible} if the associated isomorphism $\phi$ is.
\end{defin}

Let $F(X)$ be the free group on $X$.  We use the following special case of \cite[Thm 2.4]{MR1714852}:
\begin{theorem}\label{thm:StallingLinkCriterion}
	Let $F = F(S)$ be a free group with $2\leq |S|<\infty$.  Let $w$ be a cyclically reduced word in $S^{\pm}$.  Suppose $w$ contains $pq$ for each $p,q\in S^{\pm}$ with $pq$ reduced.  Then the element represented by $w$ is not contained in any proper free factor of $F$.  
\end{theorem}

We also need Lemma~\labelcref{lem:CoreDegreeCalculation} for proving Proposition~\labelcref{prop:IrredEmbedding}.

\begin{defin}
	The \emph{based core at $x$} of a connected graph $B$ is the smallest connected subgraph containing $x$ and all closed cycles of $B$. 
	Let $H\leq F$ be a finitely generated subgroup of a free group.  The \emph{core} of $H$ is the core of the associated cover of $H\leq F$.  
\end{defin}
See \cite[Constr 5.4]{MR695906} for the construction of the core of $H\leq F$ by folding.

\begin{lem}\label{lem:CoreDegreeCalculation}
	Let $F$ be a free group.  Let $A_i\in F$ be a reduced word for each $i\in I$ with $|I|<\infty$.  Let $\Gamma$ be the graph of the core associated with $\langle A_i\colon i\in I\rangle\leq F$.  Let $x$ be the basepoint of $\Gamma$.  Then $\deg x \leq 2|I|$.  
\end{lem}
\begin{proof}
	Assume $\deg x > 2|I|$.  Graph $\Gamma$ contains no vertex of degree $1$ except possibly at $x$ because each $A_\ell$ has no backtrack for $\ell\in I$.  Observe $\chi(\Gamma) \geq 1 - |I|$ since $\Gamma$ is constructed from a bouquet of $|I|$ circles and folding does not decrease $\chi(\Gamma)$.  Moreover, $\chi(\Gamma) = \sum_{v} (1 - \frac{\deg v}{2})$ since $\Gamma$ is a graph.  This leads to a contradiction since $\chi(C)\leq 1 - \frac{\deg v}{2} < 1 - |I| = \chi(\Gamma)$, where the first inequality holds since $1 - \frac{\deg v}{2}\leq 0$ for each vertex $v$. 
\end{proof}

\begin{construction}\label{construction:irredhalfclean_to_irredasc}
	Given $H$ a partial ascending fully irreducible HNN extension of a finitely generated free group $F$, we will construct a fully irreducible ascending HNN extension $G$ with $H\leq G$.  The subgroup relationship is proved in Proposition~\labelcref{prop:IrredEmbedding}.
	
	Observe $H$ can be presented as Equation~\labelcref{eq:HalfCleanHNN_GroupH}.  Recall that no $A_\ell$ has a backtrack for $\ell\in I$ since each $A_\ell$ is reduced.
	
	Let $\Gamma$ be the core of $\langle A_i \colon i\in I\rangle\leq F(a_i,b_j)$ with basepoint $x$.  Then, $\deg x \leq 2|I|$ by Lemma~\labelcref{lem:CoreDegreeCalculation}, so there are at most $2|I|$ labels in $\{a_i^{\pm}, b_j^{\pm}\colon i\in I, j\in J\}$ used on the edges emerging from $x$.  Therefore, we enumerate $2|J|$ labels in $\{a_i^{\pm}, b_j^{\pm}\colon i\in I, j\in J\}$ not used as $x_1, x_2, \dots, x_{2|J|}$.  
	
	Let $\rho\colon F(a_i, b_j, c_1, c_2)\to F(c_1, c_2)$ be the projection map.  Let $B_j, C_1, C_2\in F(a_i, b_j, c_1, c_2)$ for each $j\in J$ as follows:
	\begin{align*}
		B_j\ &=\ x_{j+|J|}\ U_j\ \beta_{j}\ x_{j}^{-1}\\
		C_1\ &=\ c_1\ V_1\ \gamma_{1}\ c_1^{-1}\\
		C_2\ &=\ c_2\ V_2\ \gamma_{2}\ c_2^{-1}
	\end{align*}
	where 
	\begin{itemize}
		\item $U_j, V_1, V_2$ are words such that $pq$ is a subword of $U_j$ and $V_k$ for each reduced $pq$ with $p, q\in \{a_i, b_j, c_1, c_2\colon i\in I, j\in J\}$, and $U_j, V_1, V_2$ do not start with $x_j^{-1}, c_1^{-1}, c_2^{-1}$ respectively;
		\item $\beta_j, \gamma_1, \gamma_2$ are words such that $W = \{\rho(B_j)\}_{j\in J}\cup \{\rho(c_1^{-1}C_1), \rho(c_2^{-1}C_2)\}$ is $C'(\frac{1}{7})$;
		\item no word in $W$ is a proper power; and
		\item $B_j, C_k$ are cyclically reduced
	\end{itemize}
	One may choose $\beta_j$ and $\gamma_k$ to be long words in $\{c_k\colon k\in K\}$ to satisfy all the conditions.
	
	We construct $G$ to be an ascending HNN extension with presentation as Equation~\labelcref{eq:HalfCleanHNN_GroupG}.
\end{construction}

\begin{prop}\label{prop:IrredEmbedding}
	With $H$ and $G$ as in Construction~\labelcref{construction:irredhalfclean_to_irredasc}, there is a group embedding from $H$ to $G$.  Moreover, $G$ is fully irreducible.
\end{prop}
\begin{proof}
	Similar to Proposition~\labelcref{prop:groupembedding}, we claim that $H=\pi_Y\to \pi_1 X = G$ is an embedding.  Indeed, observe $\overline{X}^1 = \{c_k\colon k\in K\}$ and $2$-cells in $\overline{X}$ have attaching maps labelled by $W$, thus $\overline{X}$ is $C(7)$.  Words in $W$ are all distinct and have exponent $1$.
	
	Assume $G$ is weakly reducible.  Let $F_0\leq F$ be a weakly invariant free factor of $\phi$.  Write $F = F_0*F'$.  Observe $\phi(F_0)\leq g F_0 g^{-1}$ for some $g\in F$ implies $\phi^p(F_0)$ is in some proper free factor of $F$ for some $p>0$.  Let $\Phi = \phi^p$.
	
	Since $H$ is fully irreducible, $F_0 \not\sseq F(a_i)$.  Therefore, there is $g\in F_0- F(a_i)$.  Write  $g = p_1q_1\cdots q_{n-1}p_n$ for some $p_i \in F(a_i)$ and letter $q_j \in \{b_j^\pm, c_k^\pm\}$ and $n\geq 2$ since $g\not\in F(a_i)$.  
	
	The core $\widehat{\Gamma}$ of $\langle A_i, B_j, C_k\colon i\in I, j\in J, k\in K\rangle\leq F(a_i, b_j, c_k)$ is the wedge of $\Gamma$ and bouquet of circles with paths $\{B_j, C_k\colon j\in J, k\in K\}$.  Indeed by construction, there is no folding at the basepoint when adding $B_j$ and $C_k$ for each $j\in J, k\in K$.  
	
	Therefore, $\phi(g) = \phi(p_1)\phi(q_1)\cdots \phi(q_{n-1})\phi(p_n)$ does not have folding between $\phi(p_i)$ and $\phi(q_j)$ for any $i, j$ and hence is cyclically reduced.  In particular, $U_j$ or $V_k$ is an immersed subword of $\phi(g)$.  Hence, $\phi(g)$ does not lie in any proper free factor by Theorem~\labelcref{thm:StallingLinkCriterion}, and in particular, $\phi(g)\not\in F(a_i)$.  
	
	Iterate the argument above for $p$ times to obtain $\phi^p(g)$ does not lie in any proper free factor, contradicting the free decomposition of $F$ as given by the definition.
	%
\end{proof}


\section{Relative Hyperbolicity}

We will use Osin's criterion for relative hyperbolicity \cite{MR2182268}. We first recall the definitions from Osin's paper.  Then, we focus on some types of subcomplexes of compact complexes and deduce relative hyperbolicity in that case.  

\begin{defin}
	Let $G$ be a group with a finite symmetric generating set $S$.  Let $H\leq G$ be a subgroup.  Let $\overarcmath{H}$ be an isomorphic copy of $H$.  Let $F(S\cup \overarcmath{H})$ be the free group with basis $S\cup \overarcmath{H}$.  
	We use $|w|$ to denote the length of a word $w$.  
\end{defin}

Let $\mc Z$ be the subset of $F(\overarcmath{H})$ representing the elements of $\ker(F(\overarcmath{H}) \to \overarcmath{H})$.

\begin{defin}
	Let $\mc R$ be a set of words in $S\cup \overarcmath H$ such that $G$ can be presented as $\langle S, \overarcmath{H}\mid \mc R, \mc Z\rangle$, \textit{i.e.} $F(S\cup \overarcmath{H})\to G$ has kernel $\normalclosure{\mc R, \mc Z}$.  The presentation $P = \langle S, \overarcmath{H}\mid \mc R, \mc Z\rangle$ is a \emph{finite presentation relative to $H$} if $\mc R$ is finite.  Note that $P$ is only \emph{relatively} finite.
	

	Let $\overarcmath{X}$ be the presentation complex of $P$.  Let $w$ be a word in $S\cup \overarcmath{H}$ that represents $1_G$.  The \emph{area of $w$ in $P$}, denoted $\operatorname{Area}_{_{P}}(w)$, is the minimal area of disk diagrams $D\to \overarcmath{X}$ with $\boundary_\p D = w$.  
	Suppose $P$ is a finite presentation relative to $H$.  We say $P$ has a \emph{linear isoperimetric function} if there is $K>0$ such that $\operatorname{Area}_{_{P}}(w) \leq K \cdot |w|$ for any word $w$ in $S\cup \overarcmath{H}$ representing $1_G$.  
	We say $G$ has a \emph{linear isoperimetric function relative to $H$} if there is a finite presentation $P$ relative to $H$ such that $P$ has a linear isoperimetric function.
\end{defin}

\begin{theorem}\label{thm:LinearIsoImpliesHyperbolic}
	Suppose $G$ admits a linear isoperimetric function relative to $H$.  Then $G$ is hyperbolic relative to $H$.  
\end{theorem}
\begin{proof}
	See \cite[Cor. 2.54]{MR2182268}.
\end{proof}

\begin{construction}\label{construction:FiniteRelativePresentation}
	Given a compact combinatorial $2$-complex $X$ with one $0$-cell and $\iota\colon Y\hookrightarrow X$ a subcomplex.  We will construct a finite presentation of $G = \pi_1 X$ relative to $H = \iota_*(\pi_1 Y)$.  The properties of this presentation are proved in Lemma~\labelcref{lem:CanonicalExtension}.
	
	Let $\{\bar c_i\}_i$ be the $1$-cells in $Y$.  Let $c_i = [\bar c_i]\in H$ for each $i$.  Let $\{\bar d_j\}_j$ be the $1$-cells of $X$ that are not in $Y$.  Let $d_j = [\bar d_j]\in G$ for each $j$.  
	Let $P = \langle S, \overarcmath{H}\mid \mc A, \mc B, \mc Z\rangle$, where:

	\begin{itemize}
		\item $S = \{c_i\}_i\cup \{d_j\}_j$.
		\item $\overarcmath{H}$ is a copy of $H$ where for each $h\in H$, the corresponding element is denoted $\overarcmath{h}\in \overarcmath{H}$.
		\item $\mc A = \{[\boundary_\p R]\colon R \text{ in $2$-cells of $X$}\}$. 
		\item $\mc B = \{c_i\,^{-1}{\overarcmath{c}_i}\}_i$.  
		\item $\mc Z$ is the subset of $F(\overarcmath{H})$ representing elements of $\ker(F(\overarcmath{H}) \to H)$.
	\end{itemize}
\end{construction}

\begin{rem}
	Each relator in $\mc B$ declares that $c_i$ equals its corresponding element $\overarcmath{c}_i\in \overarcmath{H}$. 
\end{rem}

\begin{lem}\label{lem:CanonicalExtension}
	With notations as in Construction~\labelcref{construction:FiniteRelativePresentation}, $P$ is a finite presentation relative to $H$.  Moreover, 
	let $\overarcmath{X}$ be the presentation complex of $P$.  There is a natural injection $X\to \overarcmath{X}$ and a retraction $r\colon \overarcmath{X}\to X$, and hence $G\cong \pi_1 \overarcmath{X}$.  
\end{lem}

\begin{proof}
	$P$ is finitely presented relative to $H$ since $X$ is compact, implying $\mc A$ and $\mc B$ are finite.   
	
	Observe that $X$ is the presentation complex of $\langle S \mid \mc A\rangle$, so there is an injection $X\to \overarcmath{X}$.  
	We now specify the retraction map $r$.  Each edge corresponding to $\overarcmath{h}\in \overarcmath{H}$ is mapped to a closed path in $X$ corresponding to $h\in H\leq G$.  Each $2$-cell of $\mc B$ maps to a backtrack $c_i^{-1}c_i$.  Each $2$-cell of $\mc Z$ maps to a disk diagram in $X$ since the boundary path corresponds to $1\in H\leq G$.  The retraction $r\colon \overarcmath{X}\to X$ implies $\pi_1 X\to \pi_1 \overarcmath{X}$ is injective. 
	
	
	Since relators in $\mc B$ equates each $c_i\in G$ to $\overarcmath{c}_i\in \overarcmath{H}$ and $\overarcmath{H}$ is generated by $\{\overarcmath{c}_i\}$, so $\pi_1 X\to \pi_1 \overarcmath{X}$ is surjective.  Therefore, $\pi_1 \overarcmath X \cong \pi_1 X = G$.
\end{proof}

\begin{defin}
	Let $\overarcmath{X}$ be the presentation complex of $\langle S, \overarcmath{H}\mid \mc A, \mc B, \mc Z\rangle$ as in Construction~\labelcref{construction:FiniteRelativePresentation}.  
	
	Let $\overarcmath D\to \overarcmath X$ be a disk diagram.  A $2$-cell $R$ in $\overarcmath D$ is a \emph{$\mc Z$-face} if $R$ maps to a $2$-cell of $\overarcmath{X}$ corresponding to an element of $\mc Z$.  The \emph{$\mc Z$-area} of $\overarcmath D$ is the number of $\mc Z$-faces in $\overarcmath D$.  We likewise define \emph{$\mc A$-face}, \emph{$\mc B$-face}, \emph{$\mc A$-area} and \emph{$\mc B$-area}.  A $1$-cell $e$ in $\overarcmath{D}$ is an $S$-edge if $e$ maps to a $1$-cell of $\overarcmath{X}$ corresponding to an element in $S$.  We likewise define $\overarcmath{H}$-edge.
	
	Since $Y$ is a subcomplex of $X$, which can be presented by $\langle S\mid \mc A\rangle$, any $\mc A$-face mapping to $Y$ is called an $\mc A_{Y}$-face and any $\mc A$-face (with interior) mapping to $X-Y$ is called an $\mc A_{X-Y}$-face.  We likewise define \emph{$\mc A_Y$-area} and \emph{$\mc A_{X-Y}$-area}.
	
	A \emph{$(\mc B + \mc Z)$-subdiagram} in $\overarcmath{D}$ is a connected component of $\left[\bcup_{\text{$\mc B$-faces }R} \widebar R\right]\cup \left[\bcup_{\text{$\mc Z$-faces }R} \widebar R\right]$.  
	
\end{defin}	

\begin{rem}
	$\boundary_\p R$ consists of only $S$-edges if $R$ is a $\mc{A}$-face; $\boundary_\p R$ consists of exactly one $S$-edge and one $\overarcmath{H}$-edge if $R$ is a $\mc{B}$-face; $\boundary_\p R$ consists of only $\overarcmath{H}$-edge if $R$ is a $\mc{Z}$-face.
\end{rem}

\begin{construction}\label{construction:Associated_Disk_Diagram}
	Given a disk diagram $\overarcmath{D}\to \overarcmath{X}$, we construct an \emph{associated disk diagram} $D\to X$.  For every $(\mc B + \mc Z)$-subdiagram $E$ in $\overarcmath{D}$, since $\boundary_\p E$ maps to labels in $S$, we may replace the subdiagram by a minimal (possibly singular) disk diagram in $X$ with $\mc A$-faces.  
\end{construction}

\begin{prop}\label{prop:sufficientrelativehyp}
	Let $Y\sseq X$ be a subcomplex of a compact $2$-complex $X$ with one $0$-cell.  
	Suppose $X/Y = \widebar X$ is combinatorially reducible and $X\to \widebar X$ has liftable cancellable pairs. 
	Suppose there is $K>0$ such that for each reduced disk diagram $D\to X$, the induced diagram $\widebar D\to \widebar X$ satisfies $\operatorname{Area}(\widebar{D})\leq K\cdot |\boundary_\p \widebar{D}|$.  
	Then $\pi_1 X$ is hyperbolic relative to $\iota_* (\pi_1 Y)$.  
\end{prop}
\begin{proof}
	By Proposition \labelcref{prop:sufficientinjection}, $\pi_1Y\to \pi_1 X$ is injective and so we identify $\pi_1Y$ with $\iota_* (\pi_1Y) \leq \pi_1 X$.  Let $G = \pi_1 X$ and $H = \iota_*(\pi_1 Y)$.   Then by Lemma~\labelcref{lem:CanonicalExtension}, $G$ admits a finite presentation $P = \langle S, \overarcmath{H}\mid \mc A, \mc B, \mc Z\rangle$ relative to $H$.  We will show that $P$ has a linear isoperimetric function, and so the result follows by Theorem~\labelcref{thm:LinearIsoImpliesHyperbolic}.  
	
	Let $\overarcmath{X}$ be the presentation complex of $P$.  Let $\sigma\to \overarcmath X$ be a nullhomotopic path.  When $\sigma$ only consists of $\overarcmath H$-edges, 
	 there is a reduced disk diagram $R$ with one $\mc Z$-face with $\boundary_\p R = \sigma$.  A linear isoperimetric inequality holds for any $K\geq 1$ in this case.  Hence, we may assume $\sigma$ has an $\mc S$-edge.  
	
	Consider a disk diagram $\overarcmath D\to \overarcmath X$ with $\boundary_\p \overarcmath D = \sigma$ such that the $(\mc A_Y, \mc A_{X-Y},\mc Z, \mc B)$-area is minimal in dictionary order.  
	
	Any $\mc A_Y$-face mapping to a $2$-cell $R$ can be replaced by a $\mc Z$-face surrounded by $\mc B$-faces.  Indeed $R\sseq Y$ so $\boundary_\p R$ can be lined by $\mc B$-faces, leaving a boundary path in $\overarcmath{H}$.  As this boundary is in $\ker(F(\overarcmath{H})\to H)$, we can fill it with a $\mc Z$-face.  Thus, by minimality, there is no $\mc A_Y$-face in $\overarcmath{D}$.  

	We now show the associated diagram $D\to X$, given by Construction~\labelcref{construction:Associated_Disk_Diagram}, is reduced.  In $\overarcmath{D}$, there is no $\mc A_Y$-face.  Therefore, after replacing each $(\mc B+\mc Z)$-subdiagram with $\mc A_Y$-faces in the construction of $D$, cancellable paris may only occur between $\mc A_{X-Y}$-faces meeting along isolated edges of the replacement diagram.  See Figure~\labelcref{fig:Impossible_Reduction}.
	However, cancellable pairs of $\mc A_{X-Y}$-faces in $D$ violate minimality of $\overarcmath{D}$.  Replacing each connected component of $\mc A_Y$-faces with a $(\mc B+\mc Z)$-subdiagram, we obtain a diagram in $\overarcmath{X}$ having a cancellable pair of $\mc A_{X-Y}$-faces, contradicting minimality of $\mc A_{X-Y}$-area.  
	
	\begin{figure}
		\centering
		\includegraphics[width = 0.55\textwidth]{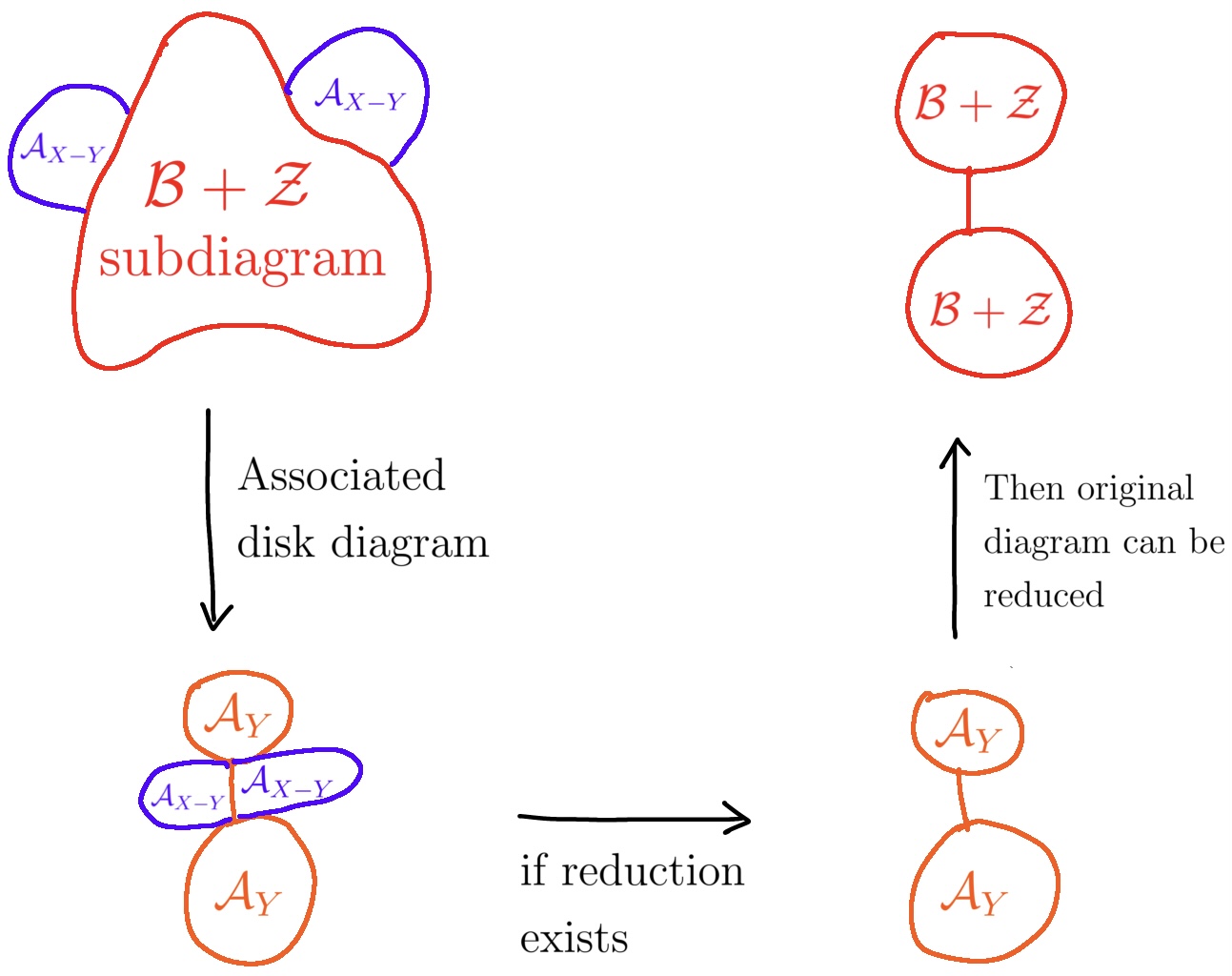}
		\caption{The configuration leading to cancellation in $D$.  The $(\mc B+\mc Z)$-subdiagram is associated with a singular disk diagram, where a pair of cancellable $\mc A_{X-Y}$-faces meet along a common edge.  Replacing each component of $\mc A_Y$-faces with a $(\mc B+\mc Z)$-subdiagram yields a diagram in $\overarcmath{X}$ with fewer $\mc A_{X-Y}$-faces, violating minimality.}
		\label{fig:Impossible_Reduction}
	\end{figure}

	Since $X\to \widebar X$ has liftable cancellable pairs, Lemma~\labelcref{lem:EquivalenceLiftableCancellablePairs} implies that $\widebar D\to \widebar X$ is reduced. The proof of Proposition~\labelcref{prop:sufficientinjection} shows that $\widebar D$ is a tree of spheres and disks.   Moreover, $\widebar D$ does not contain spheres by combinatorial reducibility of $\widebar{X}$.  
	Observe $\left|\boundary_\p \widebar{D}\right|\leq \left|\sigma\right|$ and there is no $\mc{A}_Y$-face in $\overarcmath{D}$.  The linear isoperimetric inequality for disks in $\widebar D$ implies:
	$$\# \{\text{$\mc A$-faces in $\overarcmath D$}\}
	\ =\ \# \{\text{$\mc A_{X-Y}$-faces in $\overarcmath D$}\}
	\ =\ \operatorname{Area} (\widebar {D})
	\ \leq\ K \left|\boundary_\p \widebar{D}\right|
	\ \leq\  K \left|\sigma\right|$$

	Let $M$ be an upper bound on the boundary lengths of $2$-cells of $X$.  As $\overarcmath D$ is reduced, each $\mc B$-face is either along an edge of $\boundary_\p D$ or along an $\mc A_{X-Y}$-face.  Therefore, 
	$$\# \{\text{$\mc B$-faces in $\overarcmath D$}\}
	\ \leq\ M\cdot \# \{\text{$\mc A_{X-Y}$-faces in $\overarcmath D$}\} + \left|\boundary_\p \overarcmath{D}\right|
	\ \leq\ (MK+1) \left|\sigma\right|$$

	For a $(\mc{B}+\mc{Z})$-subdiagram $E$, let $\left|E\right|_S$ be the number of $S$-edges in $\boundary E$ (possibly disconnected).  Since each $\mc B$-face has only one $S$-edge and no $\mc Z$-face has an $S$-edge, summing over all $(\mc{B}+\mc{Z})$-subdiagrams, we have:
	$\sum_i |E_i|_S 
	\ \leq\ \# \{\text{$\mc B$-faces in $\overarcmath D$}\}
	\ \leq\ (MK+1) \left|\sigma\right|$
	
	
	A $(\mc{B}+\mc{Z})$-subdiagram $E$ is homeomorphic to a genus-$0$ surface with $p(E)\geq 1$ boundary circles since $E$ is a subcomplex of $\overarcmath{D}$.  
	
	
	We claim that there is an $S$-edge in each connected component $C$ of $\boundary E$.  Indeed, assume $C$  consists only of $\overarcmath{H}$-edges.  Since no $\mc A$-face has an $\overarcmath{H}$-edge on the boundary, each $2$-cell intersecting $C$ is a $\mc B$-face or $\mc Z$-face.  Therefore, there is no $2$-cell outside $E$ since $E$ is a $(\mc B+\mc Z)$-diagram.  Hence, $C\sseq \boundary_\p D$ and since $D$ is not singular, $C = \boundary_p D = \sigma$, contradicting the assumption that $\sigma$ contains an $S$-edge.  We conclude that each $C$ contains an $S$-edge and $\left|E\right|_S\geq p(E)$.  
	
	By the construction of $\mc Z$, observe that $\mc Z$-faces can be merged if the merged cell is simply-connected.  Since $E$ is homeomorphic to a sphere with $p(E)$ boundary circles, $E$ can be merged into $p(E)$ simply-connected $\mc Z$-faces.  Let $\overarcmath{D}'$ be the disk diagram after merging $\mc Z$-faces.  We have: 
	$$\# \{\text{$\mc Z$-faces in $\overarcmath{D}'$}\}
	\ =\ \sum_{\text{$(\mc{B}+\mc{Z})$-diag }E} p(E)
	\ \leq\  \sum_{\text{$(\mc{B}+\mc{Z})$-diag }E}\left|E\right|_S
	\ \leq\ (MK+1)  \left|\sigma\right|$$
	 
	Summing up, since $\overarcmath{D}'$ is a disk diagram with $\boundary_\p \overarcmath{D}' = \sigma$, we conclude:
	\begin{align*}
		\operatorname{Area}_{_{P}}(\sigma)\ \ \leq&\ \ \#\{\text{$\mc A$-faces in $\overarcmath{D}'$}\}\ +\  \#\{\text{$\mc B$-faces in $\overarcmath{D}'$}\}\ +\  \#\{\text{$\mc Z$-faces in $\overarcmath{D}'$}\}\\
		\leq&\ \ \#\{\text{$\mc A$-faces in $\overarcmath D$}\}\ +\  \#\{\text{$\mc B$-faces in $\overarcmath D$}\}\ +\  \#\{\text{$\mc Z$-faces in $\overarcmath{D}'$}\}\\
		\leq&\ \ (2MK+K + 2) \left|\sigma\right|\qedhere
	\end{align*}

\end{proof}
\begin{cor}\label{cor:sufficientrelativehypcor}
	The statement of Proposition~\labelcref{prop:sufficientrelativehyp} holds when $X$ is a compact $2$-complex (possibly with more than one $0$-cell).
\end{cor}
\begin{proof}
	By Proposition~\labelcref{prop:sufficientinjection}, $\pi_1 Y \to \pi_1 X$ is injective and so we regard $\pi_1 Y = \iota_*(\pi_1 Y)\leq \pi_1 X$.
	
	Pick a basepoint in $Y^0$.  Choose a maximal tree $T_Y\sseq Y^1$ and extend it to $T_X\sseq X^1$.  Let $X' = X / T_X$ and $Y' = X / T_Y$.  Then $\pi_1 X' = \pi_1 X$ and $\pi_1 Y' = \pi_1 Y$.  
	
	Since $X$ is compact, $T_X$ has a finite diameter $d$.  Hence, every closed path $p'\to X'/Y'$ canonically corresponds to $p\to X/Y$  with $|p| \leq (d+1)|p'|$.  By the correspondence, the conditions of combinatorial reducibility and having liftable cancellable pairs inherit to $X'$ and $Y'$. 
	Therefore, by Proposition~\labelcref{prop:sufficientrelativehyp}, it suffices to show that there is $K>0$ such that for each reduced disk diagram $D'\to X'$, the projected diagram $\widebar{D}\to X'/Y'$ satisfies $\operatorname{Area}(\widebar{D}')\leq K\cdot \left|\boundary_\p D'\right|$.  
	
	By the correspondence between $(X,Y)$ and $(X',Y')$, let $D\to X$ correspond to $D'\to X'$ and $\widebar D\to \widebar X$ correspond to $\widebar{D}'\to X'/Y'$.  As the correspondence is canonical, $D\to X$ is reduced.  Thus, 
	$$\operatorname{Area}(\widebar{D}')
	\ =\ \operatorname{Area}(\widebar{D})
	\ \leq\ K\left|\boundary_\p D\right|
	\ \leq\ K(d+1)\left|\boundary_\p D'\right| \hfill \eqno\qedhere$$ 
\end{proof}
\begin{cor}
	$H$ is an almost malnormal quasi-isometrically embedded subgroup in $G$. 
\end{cor}

\begin{proof}
	Since $G$ is hyperbolic relative to $H$ by Proposition \labelcref{prop:sufficientrelativehyp}, so $H$ is almost malnormal 
	and quasi-isometrically embedded.  See for instance \cite[Thm 1.2]{MR2922380}.
\end{proof}

\begin{cor}
	Let $G$ and $H$ as in Constructions~\labelcref{construction:halfclean_to_ascending}~and~\labelcref{construction:irredhalfclean_to_irredasc}.  $G$ is  hyperbolic relative to $H$.  
\end{cor}
\begin{proof}
	By constructions, $\widebar X$ is a $C(7)$ presentation complex, and $X$ has no relative extra power and no duplicates relative to $Y$.  By Proposition~\labelcref{prop:C7toCombinReducible}, $\widebar X$ is combinatorially reducible.  By Lemma~\labelcref{lem:LiftableCancellable}, $X\to \widebar X$ has liftable cancellable pairs.  
	Moreover, $\widebar X$ is $C(7)$, so all reduced diagram $\widebar D\to \widebar X$ satisfies linear isoperimetric function by Proposition~\labelcref{prop:C7LinearIsoFunc}.  
	Therefore, $G$ is hyperbolic relative to $H$.  
\end{proof}

\bibliographystyle{alpha}
\bibliography{Citation}

\begin{thebibliography}{Osi06}

\bibitem[Osi06]{MR2182268}
Denis~V. Osin.
\newblock Relatively hyperbolic groups: intrinsic geometry, algebraic
  properties, and algorithmic problems.
\newblock {\em Mem. Amer. Math. Soc.}, 179(843), 2006.

\bibitem[Sta99]{MR1714852}
John~R. Stallings.
\newblock Whitehead graphs on handlebodies.
\newblock In {\em Geometric group theory down under ({C}anberra, 1996)}, pages
  317--330. de Gruyter, Berlin, 1999.

\end{thebibliography}


\begin{thebibliography}{Kam33}

\bibitem[Kam33]{VanKampen33}
Egbert R.~Van Kampen.
\newblock On the connection between the fundamental groups of some related
  spaces.
\newblock {\em American Journal of Mathematics}, 55(1):261--267, 1933.

\bibitem[Osi06]{MR2182268}
Denis~V. Osin.
\newblock Relatively hyperbolic groups: intrinsic geometry, algebraic
  properties, and algorithmic problems.
\newblock {\em Mem. Amer. Math. Soc.}, 179(843), 2006.

\bibitem[Sta83]{MR695906}
John~R. Stallings.
\newblock Topology of finite graphs.
\newblock {\em Invent. Math.}, 71(3):551--565, 1983.

\bibitem[Sta99]{MR1714852}
John~R. Stallings.
\newblock Whitehead graphs on handlebodies.
\newblock In {\em Geometric group theory down under ({C}anberra, 1996)}, pages
  317--330. de Gruyter, Berlin, 1999.

\end{thebibliography}


\begin{thebibliography}{Kam33}

\bibitem[AW23]{abdenbi_2023}
Brahim Abdenbi and Daniel~T. Wise.
\newblock Negative immersions and finite height mappings, 2023.

\bibitem[Bow12]{MR2922380}
B.~H. Bowditch.
\newblock Relatively hyperbolic groups.
\newblock {\em Internat. J. Algebra Comput.}, 22(3):1250016, 66, 2012.

\bibitem[FH99]{MR1709311}
Mark Feighn and Michael Handel.
\newblock Mapping tori of free group automorphisms are coherent.
\newblock {\em Ann. of Math. (2)}, 149(3):1061--1077, 1999.

\bibitem[Kam33]{VanKampen33}
Egbert R.~Van Kampen.
\newblock On the connection between the fundamental groups of some related
  spaces.
\newblock {\em American Journal of Mathematics}, 55(1):261--267, 1933.

\bibitem[Lyn66]{Lyndon_1966}
Roger~C. Lyndon.
\newblock On dehn’s algorithm.
\newblock {\em Mathematische Annalen}, 166(3):208–228, 1966.

\bibitem[Osi06]{MR2182268}
Denis~V. Osin.
\newblock Relatively hyperbolic groups: intrinsic geometry, algebraic
  properties, and algorithmic problems.
\newblock {\em Mem. Amer. Math. Soc.}, 179(843):1–100, 2006.

\bibitem[Sta83]{MR695906}
John~R. Stallings.
\newblock Topology of finite graphs.
\newblock {\em Invent. Math.}, 71(3):551--565, 1983.

\bibitem[Sta99]{MR1714852}
John~R. Stallings.
\newblock Whitehead graphs on handlebodies.
\newblock In {\em Geometric group theory down under ({C}anberra, 1996)}, pages
  317--330. de Gruyter, Berlin, 1999.

\end{thebibliography}
\end{document}